\documentclass[11pt]{amsart}
\usepackage[utf8]{inputenc}
\usepackage[hidelinks]{hyperref}
\usepackage[margin=1in]{geometry}
\usepackage[backend=biber,doi=false,isbn=false,url=false,eprint=false,giveninits=true]{biblatex}
\DeclareFieldFormat{pages}{#1}
\renewbibmacro{in:}{,}

\usepackage{listings}
\lstnewenvironment{code}{\noindent\minipage{\linewidth}}{\endminipage}

\usepackage{tikz}
\usepackage{pgfplots}
\pgfplotsset{compat=1.7}

\usepackage[nameinlink]{cleveref}
\newcommand{\mythm}[2]{\newtheorem{#1}[theorem]{#2} \crefalias{theorem}{#1}}

\newtheorem{theorem}{Theorem}[section]
\numberwithin{figure}{section}
\numberwithin{equation}{section}

\mythm{lemma}{Lemma}
\mythm{proposition}{Proposition}
\mythm{corollary}{Corollary}

\theoremstyle{definition}
\mythm{definition}{Definition}
\mythm{question}{Question}
\mythm{algorithm}{Algorithm}
\mythm{example}{Example}

\theoremstyle{remark}
\mythm{remark}{Remark}

\newcommand{\MM}{\mathcal{M}}
\newcommand{\NN}{\mathcal{N}}
\newcommand{\LL}{\mathcal{L}}
\newcommand{\II}{\mathcal{I}}
\newcommand{\KK}{\mathcal{K}}
\newcommand{\Z}{\mathbb{Z}}

\DeclareMathOperator{\Mon}{Mon}
\DeclareMathOperator{\maxi}{maxi}
 \title{Upper bounds for Betti numbers from constraints on the Hilbert function}
\author{Jay White}
\begin{document}
\begin{abstract}
We describe an algorithm for finding sharp upper bounds for the total Betti numbers of a saturated ideal given certain constraints on its Hilbert function.
This algorithm is implemented in the package \lstinline{MaxBettiNumbers} along with variations that produce ideals with maximal total Betti numbers.
\end{abstract}
\maketitle
\section{Introduction}
Bigatti, Hulett, and Pardue showed that the lexsegment ideal has maximum graded Betti numbers over all ideals with a given Hilbert function.
Caviglia and Murai constructed an ideal that has maximum total Betti numbers over all saturated ideals with a given Hilbert polynomial.
Both of these situations constrain the Hilbert function in some way.
We will look at a more general way of constraining the Hilbert function.

We will consider ideals in the polynomial ring $S=K[x_0,\ldots,x_{n+1}]$, but a reduction will lead to us focusing on the ring in one fewer variable, $R=K[x_0,\ldots,x_n]$, which can be embedded in $S$.
We will consider the lexicographic ordering where $x_0>\cdots>x_n$, and we say an ideal is \emph{lexsegment} if each degree contains only the lexicographically first elements in that degree.
We write $\beta_{q,j}^S(I)$ to denote the graded Betti number of an ideal, $I\subset S$, in homological degree $q$ and graded degree $j$.
Also, we write $\beta^S_q(I):=\sum_j\beta^S_{q,j}(I)$ to denote the total Betti number of $I$ in homological degree $q$.
Finally, we will use $h_{S/I}(d):=\dim_K[S/I]_d$ to denote the Hilbert function of $S/I$ in degree $d$.

Macaulay showed that any Hilbert function, $h$, has a unique lexsegment ideal, which we will denote with $L_h$.\parencite{macaulay1927} In fact, this ideal gives upper bounds for the graded Betti numbers.
This can be stated as follows:
\begin{theorem}[\cite{bigatti1993,hulett1993,pardue1996}]\label[theorem]{thm:BHP}
  Let $h$ be the Hilbert function of some quotient of $S$, and consider the following family of ideals.
  \[\II := \left\{\text{ideal }I\subset S\,\left|\,h_{S/I}(d)=h(d)\text{ for all $d$}\right.\right\}\]
  Then, the lexsegment ideal, $L_h$, with Hilbert function $h$ has maximum graded Betti numbers in $\II$.
\end{theorem}

This means that if we know the Hilbert function of an ideal, we know upper bounds for its graded Betti numbers.
Suppose, however, we only know the Hilbert polynomial.
We can do the following, similar thing for the total Betti numbers:

\begin{theorem}[\cite{cavigliamurai2013}]\label[theorem]{thm:CM}
  Let $h$ be the Hilbert function of some quotient of $S$, and consider the following family of saturated ideals.
  \[\II := \left\{\text{saturated ideal }I\subset S\,\left|\,h_{S/I}(d)=h(d)\text{ for $d\gg0$}\right.\right\}\]
  Then, there exists an ideal in $\II$ with maximum total Betti numbers in $\II$.
Furthermore, there is a construction for such an ideal.
\end{theorem}

In this paper, we will focus on the following generalization of the preceding families.
\begin{definition}
  Given a numeric function, $h:\Z\to\Z$, we define the numeric function $\Delta h:\Z\to\Z$ by $\Delta h(d) := h(d) - h(d-1)$.

  Let $G,F,g,f$ be numeric functions with $G(d)=F(d)$ for $d\geq D$, where $D$ is some integer.
  Define the following family of saturated ideals.
  \[\II(G,F,g,f) := \left\{\text{saturated ideal }I\subset S\,\left|\,\begin{array}{c}G(d)\leq h_{S/I}(d)\leq F(d)\\ g(d)\leq \Delta h_{S/I}(d)\leq f(d)\end{array}\text{ for all $d$}\right.\right\}\]
\end{definition}

\begin{remark}\label[remark]{rem:twopoints}
There are two things that we may wish to compute:
\begin{enumerate}
  \item Sharp upper bounds for the total Betti numbers of ideals in $\II(G,F,g,f)$.
  \item Ideals that have maximal total Betti numbers in $\II(G,F,g,f)$
\end{enumerate}
Notice that we use the word ``maximal'' here.
This is because some such families have no ideal with \emph{maximum} total Betti numbers.
In other words, there is not always an ideal in the family that realizes all upper bounds simultaneously.
\end{remark}
\begin{remark}
  Let $p$ be the Hilbert polynomial of some quotient of $S$.
Let $D$ be the maximum of the regularities of the (finitely many) Hilbert functions of saturated ideals that are eventually equal to $p$.
Define $G,F,g,f$ as follows:
  \begin{align*}
    G(d) &= \begin{cases}0&d< D\\p(d)&d\geq D\end{cases} & F(d) &= \begin{cases}h_S(d)&d< D\\p(d)&d\geq D\end{cases} & f(d) &= \Delta h_S(d)=h_R(d) & g(d) &= 0
  \end{align*}
\end{remark}
Then, $\II(G,F,g,f)$ is the set of saturated ideals with Hilbert polynomial $p$.
That is, it is the same family as that in \cref{thm:CM}.

\section{Reduction}\label[section]{sec:two}
We will now take a saturated ideal in $S$ and find a lexsegment ideal in $R$ with greater graded Betti numbers.
Recall that $L_h$ is the unique lexsegment ideal with Hilbert function $h$.
\begin{lemma}[\cite{cavigliamurai2013}]\label[lemma]{lem:reduction}
  Given a saturated ideal $I\subset S$, there is a lexsegment ideal $\displaystyle L:=L_{\Delta h_{S/I}}\subset R$ such that $\displaystyle h_{S/I} = h_{S/LS}$.
  Furthermore, $\displaystyle \beta^S_{q,j}(I) \leq \beta^S_{q,j}(LS)=\beta^R_{q,j}(L)$.
\end{lemma}
In other words, $L_{\Delta h_{S/I}}S$ has the same Hilbert function as $I$, and has graded Betti numbers at least as large as the graded Betti numbers of $I$.
This allows us to restrict our search to the set of Hilbert functions of quotients of $R$ rather than the set of all saturated ideals in $S$.
Naturally, we turn next to finding a formula for the total Betti numbers of a lexsegment ideal from its Hilbert function.
To do this, we need some notation.
\begin{definition}
  We will let $R_d[k]$ denote the $k$th lexicographically last monomial of $R$ in degree $d$ (where $R_d[1]=x_n^d$.) For a monomial $x_0^{a_0}\cdots x_n^{a_n}\in\Mon(R)$, we will let $\maxi(x_0^{a_0}\cdots x_n^{a_n}) := \max\{i\mid a_i> 0\}$.
  Finally, we will write
  \begin{align*}
    V_q[d,\ell] &= \sum_{k=1}^\ell \left(\binom{n+1}{q+1} - \binom{\maxi\left(R_d[k]\right)+1}{q+1}\right)
  \end{align*}
\end{definition}

We use this notation to give us the following easy-to-compute formula for total Betti numbers that is derived from the resolution in \parencite{eliahoukervaire1990}.
Recall that $\beta^S_q\left(L_{\Delta h}S\right) = \beta^R_q\left(L_{\Delta h}\right)$, so we will simply write $\beta_q\left(L_{\Delta h}\right)$ instead.
\begin{lemma}\label[lemma]{lem:formula}
  Let $\Delta h$ and $\overline{\Delta h}$ be the Hilbert functions of some quotients of $R$.
  If $\Delta h(d)=\overline{\Delta h}(d)$ for $d>D$ where $D$ is some integer, then
  \begin{align*}
    \beta_q\left(L_{\Delta h}\right) &= \beta_q\left(L_{\overline{\Delta h}}\right) + \sum_{d=0}^D V_q[d,\Delta h(d)] - \sum_{d=0}^D V_q[d,\overline{\Delta h}(d)]
  \end{align*}
\end{lemma}
\begin{remark}\label[remark]{rem:Vq}
  Let $h$ be the Hilbert function of some ideal in $\II(G,F,g,f)$.
  Because $G(d)=h(d)=F(d)$ for all $d\geq D$, we have that $\Delta G(d)=\Delta h(d)=\Delta F(d)$ for all $d>D$.
  By choosing some $\overline{\Delta h}$ that matches $\Delta h$ for $d>D$, we can use this $\overline{\Delta h}$ for all ideals in our family $\II(G,F,g,f)$.
  This means the only variable in the formula of \cref{lem:formula} is $\sum_{d=0}^D V_q[d,\Delta h(d)]$.
\end{remark}

In summary, \Cref{lem:reduction,lem:formula} allow us to maximize the total Betti numbers by maximizing $\sum V_q$ over Hilbert functions.
That is, we wish to find
\[\max\left\{\left.\sum_{d=0}^D V_q[d,\Delta h(d)]\,\right|\,h\in \LL(G,F,g,f)\right\}\]
where
\[\LL(G,F,g,f) := \left\{\text{numeric function }h\,\left|\,\begin{array}{c}G(d)\leq h(d)\leq F(d)\\ g(d)\leq \Delta h(d)\leq f(d)\\ \left[\Delta h(d)\right]_{\langle d\rangle}\leq \Delta h(d-1)\end{array}\text{ for all $d$}\right.\right\}\]
\begin{remark}\label[remark]{rem:mLB}
  Here we use the notation $a_{\langle d\rangle}$ to mean the Macaulay lower bound on the Hilbert function in degree $d-1$, given that it is $a$ in degree $d-1$.
  \footnote{
    This is in contrast with the $a^{\langle d\rangle}$ notation commonly used for the Macaulay upper bound in degree $d+1$.
    However, this is not the same notation as used in Green's proof of Macaulay's Theorem.
  }
\end{remark}
We can reformulate the two things from \cref{rem:twopoints} that we wish to compute as follows:
\begin{enumerate}
  \item Sharp upper bounds for $\sum_{d=0}^D V_q[d,\Delta h(d)]$ over $\LL(G,F,g,f)$.
  \item Functions $h\in\LL(G,F,g,f)$ with maximal $\sum_{d=0}^D V_q[d,\Delta h(d)]$.
\end{enumerate}

There is one situation where the family $\LL(G,F,g,f)$ can be simplified further.
Suppose we ``remove'' the upper bound constraint, $h_{S/I}(d)\leq F(d)$.
We can do this by making $F(d)$ as large as possible (i.e. $F(d)=h_S(d)$ or $F(d)=\infty$) when $d<D$.
\begin{proposition}[\cite{white2021}]\label[proposition]{pro:breaking}
  Consider the family $\LL(G,F,g,f)$ and suppose that $F(d)=h_{S}(d)$ for $d<D$ and $F(d)=G(d)$ for $d\geq D$ where $D$ is some integer.
  Define the following, larger, family of numeric functions that are not necessarily Hilbert functions:
  \[\NN(G,F,g,f) := \left\{\text{numeric function }h\,\left|\,\begin{array}{c}G(d)\leq h(d)\leq F(d)\\ g(d)\leq \Delta h(d)\leq f(d)\end{array}\text{ for all $d$}\right.\right\}\]
  Then, for every $h\in \NN(G,F,g,f)$, there is some $h'\in\LL(G,F,g,f)$ such that \[\sum_{d=0}^D V_q[d,\Delta h(d)]\leq \sum_{d=0}^D V_q[d,\Delta h'(d)]\]
\end{proposition}
In other words, if we don't set an upper bound for $h_{S/I}$, knowledge of the ideal structure is unneeded, and we can simply look at all numeric functions in $\NN(G,F,g,f)$.
While this family is larger, it allows us to makes simplifications that improve the search speed.

\section{Algorithm}
Throughout this section, we will consider $G,F,g,f$ to be fixed, which means that $\LL(G,F,g,f)$ and $\NN(G,F,g,f)$ are fixed as well.
We will describe an algorithm for finding a maximum for the $q$th total Betti number over the family $\II(G,F,g,f)$.
We will also briefly describe how this algorithm can be modified to find the ideals that attain maximal total Betti numbers.
These algorithms are based on the theory in \cref{sec:two}.
For more details, see \parencite{white2021}.
We will start with the situation where the Hilbert function is unbounded above (i.e.  suppose $F(d)=h_{S}(d)$ for $d<D$.)

\subsection{No upper bound for \texorpdfstring{$h_{S/I}$}{\$h\_\{S/I\}\$}}
By \cref{rem:Vq,pro:breaking}, we wish to find \[\max\left\{\left.\sum_{d=0}^D V_q[d,\Delta h(d)]\,\right|\,h\in \NN(G,F,g,f)\right\}\]
We define the following sets of tuples:
\begin{align*}
  \MM(d) &:= \left\{(\ell_0,\ldots,\ell_d)\,\left|\,
  \begin{array}{c}
    G(i)\leq \sum_{k=0}^i\ell_k\leq F(i)\\
    g(i)\leq \ell_i\leq f(i)
  \end{array}\text{ for }0\leq i \leq d\right.\right\}\\
  \MM(d, c) &:= \left\{(\ell_0,\ldots,\ell_d) \in \MM(d) \,\left|\, \sum_{i=0}^{d} \ell_i = c\right.\right\}
\end{align*}
Because $G(d)= h(d)= F(d)$ for all $d\geq D$ and $h\in \NN(G,F,g,f)$, we have that $h\mapsto (\Delta h(0),\ldots,\Delta h(D))$ is a bijection from $\NN(G,F,g,f)$ to $\MM(D)=\MM(D, G(D))$.
Thus, it is sufficient to study $\MM(D, G(D))$.

The algorithm hinges on the fact that we can write $\MM(d, c)$ as the following recursive union:
\begin{align*}
  \MM(d, c) &= \bigcup_{j=g(d)}^{f(d)}\MM(d-1,c-j)\times j
\end{align*}
where $\MM(-1,0)= \{()\}$ and $\MM(-1, c)=\emptyset$ when $c \neq 0$.
Here we use $\MM(d-1,c-j)\times j$ to mean that we append the element $j$ to each $(d-1)$-tuple to make it a $d$-tuple.

Because we are appending the value $j$ to each tuple in $\MM(d-1,c-j)$, we only need to look at the maximal tuples in each $\MM(d-1,c-j)$ rather than all of them.
More explicitly, if we define $maxV_q(d,c)$ to be the maximum of $\sum_{k=0}^d V_q[d,\ell_k]$ over $\MM(d, c)$, we get that
\begin{align*}
  maxV_q(d,c) &= \max\left\{\left. maxV_q(d - 1, c - j) + V_q[d, j] \,\right|\,  g(d)\leq j\leq f(d)\right\}
\end{align*}
Because $\NN(G,F,g,f)$ maps to $\MM(D, G(D))$, we simply need to compute $maxV_q(D, G(D))$.

The following is an algorithm that can be used to find the maximum $q$th total Betti number.
This is the algorithm used when the method \lstinline{maxBettiNumbers} has the option \lstinline{Algorithm => "Simplified"} specified.
If not specified, this algorithm will be run by default when no upper bound for $h_{S/I}$ is given.

\begin{algorithm}\label[algorithm]{alg:simple}
  Suppose $F(d)=h_{S}(d)$ for $d<D$ and $F(d)=G(d)$ for $d\geq D$, where $D$ is some integer.
  We will maximize $\sum V_q$ over $\NN(G,F,g,f)$, and thus the $q$th total Betti number over $\II(G,F,g,f)$.
  In this algorithm, $\verb`maxVDict`(d,c)$ is the value of $maxV_q(d,c)$ and $\verb`maxHFDict`(d,c)$ is the set of all tuples with that value.
  \begin{itemize}
    \item We initialize the base case by creating a dictionary, \verb`maxVDict`, containing \verb`(-1,0)=>0` and a dictionary, \verb`maxHFDict`, containing \verb`(-1,0)=>{()}`.
    \item For each value of $d$ from $0$ to $D$ do:
    \begin{itemize}
      \item Compute $V_q[d,j]$ for each $j$ from $g(d)$ to $f(d)$.
      \item For each value of $c$ from $G(d)$ to $F(d)$ do:
      \begin{itemize}
        \item Initialize \verb`maxHF` to be $\emptyset$ and $\verb`maxV`$ to be $0$.
        \item For each value of $j$ from $g(d)$ to $f(d)$ do:
        \begin{itemize}
          \item Let $c'=c-j$, $d'=d-1$.
          \item Compute $\verb`V0`=\verb`maxVDict`(d',c')+V_q[d,j]$.
          \item Compute $\verb`HF0`=\verb`maxHFDict`(d',c')\times j$.
          \item If $\verb`V0`>\verb`maxV`$ replace \verb`maxV` and \verb`maxHF` with \verb`V0` and \verb`HF0`.
          \item If $\verb`V0`=\verb`maxV`$ add the tuples in \verb`HF0` to \verb`maxHF`.
          \item If $\verb`V0`<\verb`maxV`$ do nothing.
        \end{itemize}
        \item Add the entry \verb`(d,c)=>maxV` to the dictionary \verb`maxVDict`
        \item Add the entry  \verb`(d,c)=>HF0` to the dictionary \verb`maxHFDict`
      \end{itemize}
    \end{itemize}
    \item Return the value $\verb`maxVDict`(D, G(D))$ and the set $\verb`maxHFDict`(D, G(D))$.
  \end{itemize}
\end{algorithm}

\Cref{alg:simple} will give the sharp upper bound for $\beta_q$.
We can easily modify this algorithm to run on all possible $q$ simultaneously, which would give sharp upper bounds for the total Betti numbers.
This is implemented in the package with the option \lstinline{ResultsCount => "None"}, which is the default.

There are two possible ways to find Hilbert functions with maximal total Betti numbers.
The most straightforward way is to find the functions that maximize $\sum_{q} \beta_q$, the sum of the total Betti numbers.
All functions that maximize $\sum_{q} \beta_q$ necessarily have maximal total Betti numbers.
We can further modify \cref{alg:simple} by making \verb`maxHFDict` track this sum.
The downside of this method is that not all functions with maximal total Betti numbers are discovered.
This modification is implemented in the package with the option \lstinline{ResultsCount => "One"} or \lstinline{ResultsCount => "AllMaxBettiSum"}.

The second way involves using the $(n+1)$-tuple $(V_0,\ldots,V_n)$ in the place of $V_q$.
These tuples form a poset with only some maximal elements.
The algorithm remains essentially the same, however, more care is needed when comparing \verb`V0` and \verb`maxV` as well as when adding \verb`HF0` to \verb`maxHF`.
This method, while slower, allows us to find all ideals with maximal total Betti numbers.
This modification is implemented in the package with the option \lstinline{ResultsCount => "All"}.

\subsection{Including an upper bound for \texorpdfstring{$h_{S/I}$}{\$h\_\{S/I\}\$}}

In this case, we wish to find \[\max\left\{\left.\sum_{d=0}^D V_q[d,\Delta h(d)]\,\right|\,h\in \LL(G,F,g,f)\right\}\]
We will modify the family $\MM$ to add the inequality that ensures the functions are Hilbert functions.
Recall how $a_{\langle d\rangle}$ is defined in \cref{rem:mLB}. We define
\[\KK(d, c, k) := \left\{(\ell_0,\ldots,\ell_d)\in \MM(d,c)\,\left|\, \ell_{i-1}\geq \left[\ell_i\right]_{\langle i\rangle} \text{ for $0< i\leq d$ and }\ell_d\geq k\right.\right\}\]
as before, we have that $h\mapsto (\Delta h(0),\ldots,\Delta h(D))$ is a bijection from $\LL(G,F,g,f)$ to $\KK(D, G(D), g(D))$.
In this case, the recursion is more complex:
\[\KK(d, c, k) = \bigcup_{j=\max\{g(d), k\}}^{f(d)}\KK(d-1,c-j, \left[j\right]_{\langle d\rangle})\times j\]
where $\KK(-1, 0, k)=\{()\}$ and $\KK(-1, c, k)=\emptyset$ when $c\neq 0$.
We can unravel this union further with another recursion:
\[\KK(d, c, j) = \left(\KK(d-1,c-j, \left[j\right]_{\langle d\rangle})\times j\right)\cup\KK(d, c, j + 1)\]
where $\KK(d, c, j) = \emptyset$ when $j > f(d)$ and $\KK(d, c, j) = \KK(d, c, j+1)$
when $j < g(d)$.
In the same vein as before, we have
\begin{align*}
  maxV_q(d,c,j) &= \max\left\{maxV_q\left(d - 1, c - j, \left[j\right]_{\langle d\rangle}\right) + V_q[d, j], maxV_q(d, c, j)\right\}
\end{align*}

Even though the set we are searching is smaller, this algorithm is slower because we must store the maximal ideals for each triple $(d, c, k)$, while in \cref{alg:simple} we only need to store the maximal ideals for each pair $(d, c)$.

This is the algorithm used when the method \lstinline{maxBettiNumbers} has the option \lstinline{Algorithm => "Complete"}.
If not specified, this algorithm will be run by default when an upper bound for $h_{S/I}$ is given.

\begin{algorithm}\label[algorithm]{alg:complex}
  Suppose $F(d)=G(d)$ for large $d$, with no other restrictions on $F$.
  We will maximize $\sum V_q$ over $\LL(G,F,g,f)$, and thus the $q$th total Betti number over $\II(G,F,g,f)$.
  In this algorithm, $\verb`maxVDict`(d,c,k)$ is the value of $maxV_q(d,c,k)$ and $\verb`maxHFDict`(d,c,k)$ is the set of all tuples with that value.
  \begin{itemize}
    \item We initialize the base case by creating a dictionary \verb`maxVDict` containing \verb`(-1,0,0)=>0` and a dictionary \verb`maxHFDict` containing \verb`(-1,0,0)=>{()}`.
    \item For each value of $d$ from $0$ to $D$ do:
    \begin{itemize}
      \item Compute $V_q[d,j]$ for each $j$ from $g(d)$ to $f(d)$.
      \item For each value of $c$ from $G(d)$ to $F(d)$ do:
      \begin{itemize}
        \item Initialize \verb`maxHF` to be $\emptyset$ and $\verb`maxV`$ to be $0$.
        \item For each value of $j$ from $f(d)$ to $g(d)$, counting downward, do:
        \begin{itemize}
          \item Let $c'=c-j$, $d'=d-1$, $k'=\left[j\right]_{\langle d\rangle}$.
          \item Compute $\verb`V0`=\verb`maxVDict`(d',c',k')+V_q[d,j]$.
          \item Compute $\verb`HF0`=\verb`maxHFDict`(d',c',k')\times j$.
          \item If $\verb`V0`>\verb`maxV`$ replace \verb`maxV` and \verb`maxHF` with \verb`V0` and \verb`HF0`.
          \item If $\verb`V0`=\verb`maxV`$ add the tuples in \verb`HF0` to \verb`maxHF`.
          \item If $\verb`V0`<\verb`maxV`$ do nothing.
          \item Add the entry \verb`(d,c,j)=>maxV` to the dictionary \verb`maxVDict`
          \item Add the entry  \verb`(d,c,j)=>HF0` to the dictionary \verb`maxHFDict`
        \end{itemize}
      \end{itemize}
    \end{itemize}
    \item Return the value $\verb`maxVDict`(D, G(D), g(D))$ and the set $\verb`maxHFDict`(D, G(D), g(D))$.
  \end{itemize}
\end{algorithm}
As before, there are several variations we can have of this algorithm.
The same variations listed after \cref{alg:simple} are included in the package.

\section{Examples}
We will consider an example where $S$ is the polynomial ring in $5$ variables ($n=3$).
This example has only maximal total Betti Numbers, and not maximum total Betti numbers.
Also, \cref{alg:simple,alg:complex} give different results.
Both of these situations are somewhat unusual, but give an illuminating example.
For this example, we will choose the following constraints:
\begin{align*}
  h_{S/I}(6) &= 41 & \Delta h_{S/I}(3) &\geq 8 & \Delta h_{S/I}(5) &\geq 5\\
  h_{S/I}(d) &= 49\text{ for }d\gg0 & \Delta h_{S/I}(4) &\geq 8 & \Delta h_{S/I}(6) &\geq 5
\end{align*}
Since we will be using these constraints in several examples, we will first define a few variables to reduce repetition.
By leaving certain degrees empty, we request the default, trivial constraints on those degrees.

\begin{code}
i1 : N = 5;
i2 : g = HilbertDifferenceLowerBound => {,,,8,8,5,5};
i3 : G = HilbertFunctionLowerBound => {,,,,,,41};
i4 : F = HilbertFunctionUpperBound => {,,,,,,41};
i5 : p = HilbertPolynomial => 49;
\end{code}

We can use \lstinline{maxBettiNumbers} to find that $(23, 54, 47, 14)$ is the upper bound for the total Betti numbers of all saturated ideals with these constraints.
Additionally, the maximum for the sum of the total Betti numbers is $137$.
Because $23+54+47+14=138$, there is no single ideal with total Betti numbers of $(23, 54, 47, 14)$.
This is indicated by the value of \lstinline{isRealizable}.

\begin{code}
i6 : maxBettiNumbers(N,p,g,G,F)
o6 = MaxBetti{BettiUpperBound => {23, 54, 47, 14}}
              isRealizable => false
              MaximumBettiSum => 137
\end{code}

If we want the Hilbert function of an ideal with maximal total Betti numbers, we can pass \lstinline{ResultsCount=>"One"} as an option.
This gives an ideal that realizes the maximum for the sum of the total Betti numbers.

\begin{code}
i7 : maxBettiNumbers(N,p,g,G,F, ResultsCount=>"One")
o7 = MaxBetti{BettiUpperBound => {23, 54, 47, 14}                         }
              HilbertFunctions => {{1, 5, 11, 21, 30, 36, 41, 46, 49, 49}}
              isRealizable => false
              MaximumBettiSum => 137
\end{code}

If we want the Hilbert function of all ideals that realize the maximum sum of the total Betti numbers, we can pass \lstinline{ResultsCount=>"AllMaxBettiSum"} as an option.

\begin{code}
i8 : maxBettiNumbers(N,p,g,G,F, ResultsCount=>"AllMaxBettiSum")
o8 = MaxBetti{BettiUpperBound => {23, 54, 47, 14}                                                 }
              HilbertFunctions => {{1, 5, 11, 21, 30, 36, 41, 42, 43, 44, 45, 46, 47, 48, 49, 49}}
                                  {{1, 5, 11, 21, 30, 36, 41, 43, 44, 45, 46, 47, 48, 49, 49}    }
                                   --------------------------14 others--------------------------
                                  {{1, 5, 11, 21, 30, 36, 41, 45, 49, 49}                        }
                                  {{1, 5, 11, 21, 30, 36, 41, 46, 49, 49}                        }
              isRealizable => false
              MaximumBettiSum => 137
\end{code}

Finally, if we want the Hilbert function of all ideals that have maximal total Betti numbers, we can pass \lstinline{ResultsCount=>"All"} as an option.
In this case, the
maximal total Betti numbers are $(23,54,45,13)$ and $(22,54,47,14)$.

\begin{code}
i9 : maxBettiNumbers(N,p,g,G,F, ResultsCount=>"All")
o9 = MaxBetti{BettiUpperBound => {23, 54, 47, 14}                                                 }
              HilbertFunctions => {{1, 5, 15, 23, 31, 36, 41, 42, 43, 44, 45, 46, 47, 48, 49, 49}}
                                  {{1, 5, 11, 21, 30, 36, 41, 42, 43, 44, 45, 46, 47, 48, 49, 49}}
                                   --------------------------32 others--------------------------
                                  {{1, 5, 15, 23, 31, 36, 41, 46, 49, 49}                        }
                                  {{1, 5, 11, 21, 30, 36, 41, 46, 49, 49}                        }
              isRealizable => false
              MaximalBettiNumbers => {{23, 54, 45, 13}}
                                     {{22, 54, 47, 14}}
              MaximumBettiSum => 137
\end{code}

The package \lstinline{MaxBettiNumbers} also provides methods that can be useful for dealing
with Hilbert functions, notably \lstinline{almostLexBetti} and \lstinline{almostLexIdeal}, which
give, respectively, the Betti table and saturated ideal with largest graded Betti numbers associated with a given Hilbert function.
These match, as expected, with \lstinline{MaximalBettiNumbers} above.

\begin{code}
i10 : almostLexBetti(N, last o9.HilbertFunctions)
             0  1  2  3  4
o10 = total: 1 22 54 47 14
          0: 1  .  .  .  .
          1: .  4  6  4  1
          2: .  .  .  .  .
          3: .  6 14 11  3
          4: .  5 14 13  4
          5: .  2  5  4  1
          6: .  .  .  .  .
          7: .  2  6  6  2
          8: .  3  9  9  3
i11 : almostLexIdeal(QQ[x_1..x_N], last o9.HilbertFunctions)
              2                     4   3     3     2 2   2       2 2     4 
o11 = ideal (x , x x , x x , x x , x , x x , x x , x x , x x x , x x , x x ,
              1   1 2   1 3   1 4   2   2 3   2 4   2 3   2 3 4   2 4   2 3 
         3       2 2       3     4   6   5     4 4   3 5   2 7     8   9
      x x x , x x x , x x x , x x , x , x x , x x , x x , x x , x x , x )
       2 3 4   2 3 4   2 3 4   2 4   3   3 4   3 4   3 4   3 4   3 4   4
\end{code}

Because we are setting an upper bound of $41\geq h_{S/I}(6)$, we cannot use the algorithm described in \cref{alg:simple}, but must use the algorithm described in \cref{alg:complex}.
The method \lstinline{maxBettiNumbers} automatically chooses the appropriate algorithm, but we can force it to use a different one if we wish.
In this case, if we specify \lstinline{Algorithm=>"Simplified"}, we get an upper bound that is, necessarily, larger.
Additionally, we are given a Hilbert function that appears to be valid, but is not the Hilbert function of any saturated ideal.

\begin{code}
i12 : maxBettiNumbers(N,p,g,G,F, Algorithm=>"Simplified", ResultsCount=>"One")
o12 = MaxBetti{BettiUpperBound => {24, 57, 50, 15}                     }
               HilbertFunctions => {{1, 5, 11, 21, 30, 36, 41, 49, 49}}
               isRealizable => false
               MaximumBettiSum => 145
\end{code}

We can compare the speed of \cref{alg:simple} and \cref{alg:complex} with an example of fixing the Hilbert polynomial to be $3d^2-6d+175$ in a ring with $6$ variables.
Because there is no upper bound for $h_{S/I}$, both algorithms give valid results, which we do not display.

\begin{code}
i13 : N = 6;
i14 : QQ[i]; p = HilbertPolynomial => 3*i^2-6*i+175;
i16 : time maxBettiNumbers(N, p, Algorithm=>"Simplified", ResultsCount=>"None");
     -- used 5.82732 seconds
i17 : time maxBettiNumbers(N, p, Algorithm=>"Simplified", ResultsCount=>"All");
     -- used 7.40666 seconds
i18 : time maxBettiNumbers(N, p, Algorithm=>"Complete", ResultsCount=>"None");
     -- used 21.4503 seconds
i19 : time maxBettiNumbers(N, p, Algorithm=>"Complete", ResultsCount=>"All");
     -- used 28.4778 seconds
\end{code}

One alternative to the algorithms in this paper is a function such as \lstinline{stronglyStableIdeals} in the package \lstinline{StronglyStableIdeals}.
This method generates all strongly stable ideals with a specified Hilbert polynomial.
It is known that strongly stable ideals have the largest graded Betti numbers.
We could then select the ideals that satisfy the constraints, compute the total Betti numbers, and then find those with maximal total Betti numbers.
We will use \lstinline{stronglyStableIdeals} as a benchmark for comparing speeds with existing methods.
Although \lstinline{stronglyStableIdeals} only generates the ideals, we are not adding to its timing the additional time it would take to compute total Betti numbers and find maximums.
As a result the actual speed differences are even more pronounced.
We will time examples with a constant Hilbert polynomial.

\begin{code}
i20 : loadPackage "StronglyStableIdeals";
i21 : benchmark("maxBettiNumbers(5, HilbertPolynomial => 25)")
o21 = .0894562386249997
i22 : benchmark("stronglyStableIdeals(25, 5)")
o22 = 61.856958304
\end{code}

The following is a graph of the timing for small constant Hilbert polynomials.
Notice that it a log plot, and that \lstinline{stronglyStableIdeals} quickly reaches over 1000 seconds.
\begin{center}
\begin{tikzpicture}
\begin{axis}[
title=Polynomial vs Runtime for Different Algorithms,
ymode=log,
xlabel=Hilbert Polynomial,
ylabel=Runtime (seconds),
cycle list name=mark list,
ymax=1000,
ymin=0.002,
xmin=-2,
xmax=102,
width=0.99\textwidth,
height=0.5\textwidth
]
\addplot+[mark size=0.5pt,mark=+] table {
0 .00888917
1 .00920395
2 .0103165
3 .0111199
4 .0123516
5 .0134797
6 .0147627
7 .0163151
8 .0182858
9 .0203721
10 .0222417
11 .024596
12 .0273881
13 .0298749
14 .0331405
15 .0363815
16 .0395687
17 .0428068
18 .0465802
19 .0504264
20 .0549085
21 .0595622
22 .0636117
23 .0687192
24 .0732659
25 .078373
26 .0845019
27 .0897713
28 .0960005
29 .101652
30 .108383
31 .115804
32 .120993
33 .127969
34 .137225
35 .144347
36 .150877
37 .159227
38 .170027
39 .174464
40 .188171
41 .194686
42 .204767
43 .212199
44 .223636
45 .233092
46 .247319
47 .255479
48 .261461
49 .275814
50 .288269
51 .297538
52 .310882
53 .319917
54 .334412
55 .34365
56 .360538
57 .373788
58 .379799
59 .398144
60 .412452
61 .430293
62 .444451
63 .457631
64 .468709
65 .488527
66 .498112
67 .523378
68 .536941
69 .552648
70 .556136
71 .579247
72 .606404
73 .616413
74 .635165
75 .643716
76 .671194
77 .698251
78 .707129
79 .724008
80 .750678
81 .759451
82 .787598
83 .797579
84 .827648
85 .842981
86 .87188
87 .885343
88 .909235
89 .929166
90 .944998
91 .978678
92 1.00013
93 1.01877
94 1.0531
95 1.03459
96 1.08842
97 1.10997
98 1.12008
99 1.14634
100 1.17162
};
\addlegendentry{Algorithm 3.1}
\addplot+[mark size=0.75pt,mark=*] table {
0 .00876827
1 .00925316
2 .0101253
3 .0109037
4 .0119534
5 .0130011
6 .0145054
7 .0161449
8 .0177849
9 .0200174
10 .0221154
11 .025222
12 .0281253
13 .0314201
14 .0351755
15 .0391211
16 .044862
17 .0482551
18 .0530953
19 .0588087
20 .0648489
21 .0704615
22 .0785537
23 .086873
24 .0938721
25 .102029
26 .112268
27 .121131
28 .132679
29 .143802
30 .153622
31 .166665
32 .180585
33 .193147
34 .209144
35 .223184
36 .240066
37 .257559
38 .27382
39 .288979
40 .311107
41 .327037
42 .353221
43 .368889
44 .390682
45 .423002
46 .442066
47 .464543
48 .491553
49 .523551
50 .547263
51 .57074
52 .600179
53 .626775
54 .666153
55 .704896
56 .734302
57 .770134
58 .812256
59 .851618
60 .896631
61 .924832
62 .974537
63 1.0136
64 1.04018
65 1.08401
66 1.12085
67 1.2001
68 1.23796
69 1.2937
70 1.33462
71 1.40439
72 1.44549
73 1.48333
74 1.57154
75 1.63619
76 1.67578
77 1.73588
78 1.78423
79 1.85601
80 1.93871
81 2.03635
82 2.07789
83 2.14066
84 2.21471
85 2.30677
86 2.37115
87 2.43094
88 2.50424
89 2.62444
90 2.71185
91 2.7856
92 2.86441
93 2.94246
94 3.07105
95 3.14318
96 3.23666
97 3.32036
98 3.4507
99 3.52319
100 3.64643
};
\addlegendentry{Algorithm 3.2}
\addplot+[mark size=1.5pt,mark=square] table {
1 .00412034
2 .00494015
3 .0076638
4 .0116511
5 .0190108
6 .0294711
7 .0486676
8 .076634
9 .121402
10 .184871
11 .282395
12 .426749
13 .639348
14 .963555
15 1.40449
16 2.04883
17 2.99862
18 4.36594
19 6.29752
20 9.0339
21 12.8946
22 18.4646
23 25.972
24 38.2488
25 56.4104
26 76.7054
27 105.703
28 140.581
29 191.014
30 258.073
31 350.746
32	517.541
33	709.804
34	972.779
35	1363.46
};
\addlegendentry{stronglyStableIdeals}
\end{axis}
\end{tikzpicture}
\end{center}
 
\printbibliography

\end{document}